\documentclass[journal ]{new-aiaa}
\usepackage[utf8]{inputenc}
\usepackage{textcomp}

\usepackage{graphicx}
\usepackage{amsmath}
\usepackage{varwidth}
\usepackage{comment}

\newtheorem{definition}{Definition}

\newtheorem{proposition}{Proposition}

\usepackage[version=4]{mhchem}
\usepackage{siunitx}
\usepackage{longtable,tabularx}
\title{Algebraic Lyapunov Functions for Homogeneous Dynamic Systems}

\author{Hassan Abdelraouf\footnote{PhD candidate, Aerospace Engineering, University of Illinois at Urbana-Champaign, Illinois, USA, hassana4@illinois.edu.}}
\author{Eric Feron\footnote{Professor, Department of  Electrical Engineering, King Abdullah University of Science and Technology, Thuwal, KSA,  eric.feron@kaust.edu.sa, AIAA Fellow.}}
\author{Jeff Shamma \footnote{Professor, Industrial and Enterprise Systems Engineering, Univresity of Illinois at Urbana-Champaign, Illinois, USA,  jshamma@illinois.edu}}
\affil{}

\begin{document}

\maketitle

\begin{abstract}
A method for constructing homogeneous Lyapunov functions of degree 1 from polynomial invariant sets is presented for linear time varying systems, homogeneous dynamic systems and the class of nonlinear systems that can be represented as such. The method allows the development of Lyapunov functions that are not necessarily symmetric about the origin, unlike the polynomial, homogeneous Lyapunov functions discussed in the prior literature. The work is illustrated by very simple examples.

\end{abstract}

\section{Introduction}
\label{sec:introduction}
\lettrine{T}{he} analysis of linear time-invariant (LTI) systems has been closely tied to quadratic Lyapunov functions ever since it was shown that LTI system stability is equivalent to the existence of such a Lyapunov function \cite{khalil2002nonlinear}. The larger class of time-varying and switching linear systems, which can encompass nonlinear systems that contain sector-bounded nonlinearities as well as systems that admit global linearizations \cite{boyd1994linear}, has also considerably benefited from the search for quadratic Lyapunov functions that prove the stability of these systems. Stability for these systems does not necessarily imply the existence of a quadratic Lyapunov function, but it is equivalent to the existence of a homogeneous polynomial Lyapunov function \cite[Theorem 3.7] {ahmadi2017sum} \cite{mason2006common}. Even more recently, there has been the development of methods allowing the computation of such homogeneous polynomial Lyapunov functions via the computation of quadratic Lyapunov functions via the analysis of a sequence of increasingly large systems, lifted from the original system~\cite{abate2020lyapunov}. For linear systems alone, the computation of such higher-order Lyapunov functions can considerably improve the rigorous analysis of their pointwise-in-time properties, such as peak impulse response, overshoot, and other quantities of interest, far beyond the capabilities offered by quadratic Lyapunov functions.

In this report, we introduce a new class of Lyapunov functions that is motivated by recent progress with the computation of polynomial invariant sets that are not necessarily symmetric about the origin, using techniques based on the relaxation of the corresponding problems to convex, semi-definite programs under the generally accepted "Sum-of-Squares" label. These Lyapunov functions are implicitly defined as the unique non-negative solution to a polynomial equation in a single variable parameterized by the system's state variables. The name "Algebraic Lyapunov function" follows from its definition.

\section{Preliminaries}
Here we consider a  dynamic system of the kind
\begin{equation}
\displaystyle \frac{d}{dt} x = f(x) , \;\; f(0)=0, \;\; x(0) = x_0,
\label{system}
\end{equation}
where $x\in \mathbb{R}^n$ is the state and $f: \mathbb{R}^n \to \mathbb{R}^n$ is a homogeneous function. 

\begin{definition}
For a function $f: \mathbb{R}^n \to \mathbb{R}^n$, if $ \forall \lambda $ and $ \forall x $ is homogeneous of degree $\nu$ if for every $\lambda$ and for every $x\in \mathbb{R}^n$ 
$f(\lambda x) =\lambda^{\nu+1} x$.
\end{definition}
\noindent For example, the linear system 
\begin{equation}
\displaystyle \frac{d}{dt} x = A x ,  \;\; x(0) = x_0,
\label{linear system}
\end{equation}
is homogeneous of degree of zero. 
It is well-known that system~(\ref{linear system}) is stable, that is, the eigenvalues of $A$
all have negative real part if and only if there exists a matrix $P$, symmetric and positive-definite, such that the matrix $A^TP + PA$ is negative definite. Then, $V(x) = x^TPx$ is a quadratic Lyapunov function for the system. The powerful equivalence between linear system stability and the existence of quadratic Lyapunov functions has led to the creation of a large number of extensions for the analysis of time-varying, and switching linear systems, on the one hand, and the analysis of input-output properties, such as passivity, non-expansivity and robustness properties of linear systems on the other hand~\cite{boyd1994linear}. However, quadratic Lyapunov functions do not constitute a panacea for the analysis of many other properties of linear systems, as shown in other recent publications~\cite{klett2022towards}. Notably homogeneous polynomial Lyapunov functions allow significant refinements of the static, that is, non-simulation based, analysis of pointwise-in-time properties of~\ref{linear system}. Motivated by this past observation, we are always seeking the formulation of ever-more flexible Lyapunov functions that aim at showing important properties of linear systems.

\section{Invariant polynomial sets}
One of the recent and valuable efforts has been the computation not only of Lyapunov functions, but also of invariant sets. Invariant sets are valuable to bound, for example, the behavior of linear systems subject to persistent, bounded excitation~\cite{abdelraouf2022computing}, as well as the behavior of chosen nonlinear systems, such as the vanderpol oscillator ~\cite{henrion2013convex}. Owing to the powerful sum-of-square relaxations associated with algebraic or semi-algebraic optimization problems, polynomial expressions of these invariant sets constitute a strong current corpus of ongoing research.

With the system~\ref{system}, we might be interested in bounding, for example, the trajectory followed by the state $x$, given the initial condition $x_0$. Consider a candidate, bounded level set $P(x) = 0$, where $P$ is a polynomial in the entries of $x$ of given degree. Assuming $P(x) < 0$ when $x$ is "inside" the level set and $P(x)>0$ when $x$ is "outside" the level set, the level set is said invariant if the trajectories of the system~\ref{system} cross the set "from outside to inside", which can be written
\begin{equation}
\frac{d}{dt}P(x(t))<0 \mbox{ whenever } P(x) = 0.
\label{crit_cond}
\end{equation}
To this condition, one must add that $P(x_0)\leq 0$ so as to ensure
that the invariant set contains the initial condition. Various techniques exist to find such invariant sets. One is to "relax" the condition~(\ref{crit_cond}) and replace it by 
\[
U_1(x)\displaystyle \frac{\partial P}{\partial x}f(x) + U_2(x)P(x) <0, \;\;\forall x \]
where $U_1$ is a positive polynomial of given degree and $U_2$ is an arbitrary polynomial of given degree. Adding the constraints that $U_1$ and $-\left(U_1(x)\displaystyle \frac{\partial P}{\partial x}f(x) + U_2(x)P(x)\right)$ be "sum of square" polynomials greatly simplify looking for appropriate $P$, $U_1$, and $U_2$ and therefore solve the problem. 

\section{Contribution: Algebraic Lyapunov functions}
Assuming an invariant polynomial set characterized by the polynomial equation $P(x) = 0$ has been computed, the question arises about extracting a Lyapunov function from this set. For that, we rewrite $P(x)$ as the sum of its monomials of given degree, that is
\[
P(x) = \displaystyle \sum_{i=0}^p M_i(x),
\]
where the polynomials $M_i(x)$ are homogeneous of degree $i$. Next, given the real variable $\tau$ we introduce the polynomial 
\[
\tilde{P}(x,\tau) = \displaystyle \sum_{i=0}^p M_i(x)\tau^{p-i}
\]
in such a way that $\tilde{P}(x,\tau)$ is homogeneous of degree $p$ in the variable $(x,\tau)$. Note that $\tilde{P}(x, 1) = P(x) = 0$ when $x$ is in the invariant set. Consider now the sets $\tilde{P}(x, \tau) = 0$ as $\tau$  ranges from zero to $\infty$. It is easy to see that $\tilde{P}(x, \tau) = 0$ is a "scaled up" or "scaled down" version of $P(x)=0$ depending on whether $\tau >1$ or $\tau<1$. For example, when $\tau = 2$, $\tilde{P}(2x,2)=0$ if $P(x)=0$.
Conversely, let $\tau(x)$ be a positive solution to the polynomial equation in one variable $\tilde{P}(x,\tau)=0$ with $x \neq 0$ given. We claim $\tau(x)$ is a homogeneous Lyapunov function of degree $1$ for the system~(\ref{system}). For this, we may need a few technical conditions to be met. 

{\em Unicity condition:} For any given $x$, we assume that there exists a unique solution to $\tilde{P}(x,\tau)=0$. It can be easily seen that an equivalent condition is that the set $P(x)$ be star-convex about the origin \cite{shamma2003existence}. Once this is assured, we show that $\tau(x)$ is indeed a Lyapunov function. 

\begin{proposition}
The function $\tau(x)=\{c>0:\tilde{P}(x,c)=0\}$ is a Lyapunov function for the system \ref{system}.
\end{proposition}
\textit{Proof:} First, the function $\tau(x)$ can be written in terms of $P(x)$ as $\tau(x)=\{c>0:{P}(x/c)=0\}$. Along the system trajectory, $x$ and $c$ are changing with time, but $P(x/c)=0$. By taking the derivative of $P(x/c)$ w.r.t. time along the system trajectory, we get  
\begin{equation}
\left(\nabla_{(x/c)} P\right)^T \left(\frac{\dot{x}}{c}-\frac{x\dot{c}}{c^2}\right)=0\\
\label{derivative1}
\end{equation}
then, 
\begin{equation}
\left(\nabla_{(x/c)} P\right)^T \left(\frac{\dot{x}}{c}\right)=\left(\nabla_{(x/c)} P\right)^T \left(\frac{x}{c}\right) \frac{\dot{c}}{c}
\label{derivative2}
\end{equation}
Let $y=x/c$, then the left hand side in \ref{derivative2} is the time derivative of $P(.)$ at the boundary of the invariant set, so it is a negative term. Therefore, 
\begin{equation}
\left(\nabla_{y} P\right)^T y \frac{\dot{c}}{c}<0
\label{RHS}
\end{equation}
since, we take the gradient of $P(x/c)$ w.r.t. $(x/c)$, we can ignore $M_0$. Additionally, the invariant set contains the origin, so $P(x/c)-M_0>0$ and can be written as $z^T Q z $ where $Q$ is a positive definite matrix and $z$ is the vector of monomials up to degree $p/2$ \cite{ahmadi2018sum}. without loss of generality, we can assume in the proof a second order system (i.e. $x\in \mathbb{R}^2$). To write $z$ in the general form, we define $v_n \in \mathbb{R}^{n+1}$ as the vector of  monomials of degree $n$, so 
$v_n = (\begin{matrix}
y_1^n & y_1^{n-1} y_2 &\dots & y_2^n\end{matrix})^T$
Hence, $z \in \mathbb{R}^{(n/2+1)(n/2+2)/2 -1}$ can be writen as $
z = (\begin{matrix}
v_{n/2}^T & v_{n/2-1}^T&\dots &v_1^T
\end{matrix})^T$. For example, if the $P(y)$ is a polynomial of  degree 4, then $z$ is the vector of monomials up to degree $2$, so $z\in \mathbb{R}^5$ and can be written as 
$z= (\begin{matrix}
y_1^2 & y_1 y_2 & y_2^2 & y_1 &y_2
\end{matrix})^T$.

Hence, the term $\left(\nabla_{y} P\right)^T y$ in the right hand side in \ref{derivative2} can be simplified as 
\begin{equation}
\begin{aligned}
\left(\nabla_{y} P\right)^T y &= y^T \nabla_{y} P \\
&=y^T \left(\frac{dz}{dy}\right)^T \frac{dV}{dz} \\
&= 2 y^T \left(\frac{dz}{dy}\right)^T Q z \\
&= 2 z^T TQ z \qquad\text{ where } T = \text{diag}\left( n/2I_{n/2+1},(n/2-1)I_{n/2},\dots,I_2\right).
\end{aligned}
\end{equation}
where $I_i \in \mathbb{R}^{i\times i}$ is the identity matrix. 

Since $Q$ is a positive definite matrix and $T$ is also a positive definite matrix because it is a diagonal matrix with positive integers on its diagonal. Therefore, the term $\left(\nabla_{y} P\right)^T y >0$ for any nonzero $y$. From \ref{RHS}, we conclude that $\dot{c}/c<0$ along the system trajectory. Since $c>0$ form the definition of $\tau(x)$, $\dot{c}<0$ which proves that $\tau(x)$ is a Lyapunov function for the system \ref{system}.

\section{Example}
Considering the following linear system 
\begin{equation}
\frac{d}{dt}\left[\begin{array}{c}x_1 \\ x_2\end{array}\right] = \left[\begin{array}{cc} -1 & 0 \\
0 & -1\end{array} \right]\left[\begin{array}{c}x_1 \\ x_2\end{array}\right],
\label{identity}
\end{equation}
the set of all $(x_1,x_2)$ satisfying
\[
P(x_1,x_2)= (x_1 -1)^2 + (x_2+1)^2 -4 =0
\]
is invariant (trivial). The polynomial $\tilde{P}$ then satisfies
\[
\tilde{P}(x_1,x_2,\tau) = -2\tau^2 + 2\tau(-x_1+x_2) + (x_1^2+x_2^2)
\]
Given $x_1$ and $x_2$, solving $\tilde{P}(x_1,x_2,\tau)=0$ yields 
\[
\tau = \left(x_2-x_1 \pm \sqrt{(-x_1+x_2)^2+2(x_1^2+x_2^2)}\right)/2,
\]
whose only positive value is 
\begin{equation}
\tau (x_1,x_2) = \left(x_2-x_1 + \sqrt{(-x_1+x_2)^2+2(x_1^2+x_2^2)}\right)/2.
\label{lyap}
\end{equation}
and it is a Lyapunov function for System~(\ref{identity}).
A graphical view of the system, together with level sets of the algebraic Lyapunov function for $\tau \in \left\{0.25, 0.5,1\right\}$

\newpage
\begin{figure}[h!]
\includegraphics[width=0.6\linewidth]{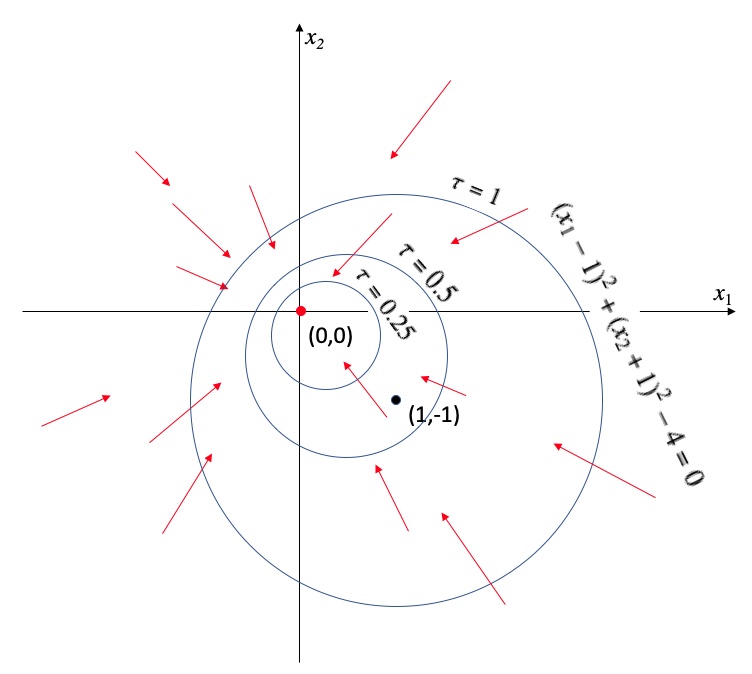}
\centering
\caption{Algebraic Lyapunov function~(\ref{lyap}) contours and System~(\ref{identity}).}
\label{fig:Vee}
\end{figure}

\section{Conclusion}
This report introduces the concept of algebraic Lyapunov function, that derives directly from recent progress in the computation of invariant sets for dynamical systems, most notably linear time-invariant or time-varying or switching systems and nonlinear systems that may be represented as such. The Lyapunov function is extracted as the unique positive root of a polynomial that is a function of the state variables, thus the name "algebraic Lyapunov function".

\section*{Funding Sources}

This work is supported by the KAUST baseline budget.

\bibliography{sample}

\end{document}